\documentclass[draft]{amsart}
\usepackage{amsmath,amsthm,amssymb}

\newtheorem{thm}{Theorem}[section]

\theoremstyle{definition}

\theoremstyle{remark}

\numberwithin{equation}{section}

\begin{document}

\title{A Torelli--like Theorem for Smooth Plane Curves}


%    author one information
\author{James S.~Wolper}
\address{Dept.~of Mathematics, Idaho State University, 921 S.~8th
Ave., Mail Stop 8085, Pocatello, ID~~83209~USA}
\email{wolpjame@isu.edu}
%\thanks{}

%    \subjclass is required.

\date{17 December, 2012}

%\dedicatory{}

%    "Communicated by" -- provide editor's name; required.
%\commby{}

%    Abstract is required.
\begin{abstract}
The Information-Theoretic Schottky Problem treats the period matrix of a compact Riemann Surface as a compressible signal. In this case, the period matrix of a smooth plane curve is characterized by only 4 of its columns, a significant compression.
\end{abstract}

\maketitle
%%%%%%%%%%%%%%%%%%%%%%%%%%%%%
\renewcommand{\sectionmark}[1]{}

%%%%%%%%%%%%%%%%%%%%%%%%%%%%%%%%%%%%%%%%%%%
\section{Introduction}

Begin by fixing notation; consult \cite{GH} as a general reference.

Let $X$ be a compact Riemann Surface 
of genus $g > 1$; equivalently,
$X$ is a nonsingular complex algebraic curve.  Choose a basis $\omega_1,
\ldots, \omega_g$ for the space $H^{1,0}(X)$ of holomorphic differentials
on $X$, and a symplectic basis $\alpha_1, \ldots, \alpha_g$, $\beta_ 1,
\ldots, \beta_g$ for the singular homology $H_1(X, {\bf Z})$, normalized
so $\int_{\alpha_i} \omega_j = \delta_{ij}$, the Dirac delta.  The 
matrix $\Omega_{ij} := \int_{\beta_i} \omega_j$ is the {\sl period
matrix\/}; Riemann proved that it is symmetric with positive definite imaginary
part.  The torus ${\mathbb C}^g / [I | \Omega]$ is the {\sl Jacobian\/} of $X$.
Torelli's Theorem asserts that the Jacobian determines all of the
properties of $X$. In practice deciding which properties apply
is seldom successful (but see \cite{W07}).

The period matrix is symmetric with positive-definite imaginary part,
and the space of such matrices forms the {\sl Siegel upper half-space\/} ${H}_g$.
Its dimension is $g(g+1)/2$, while the dimension of the moduli space of 
curves of degree $g$ has dimension $3g-3$.  Distinguishing the period
matrices from arbitrary elements of ${H}_g$ is the {\sl Schottky
Problem\/}.  See \cite{G} for details on the problem and some of its
previous solutions.

Now, recast the problem in terms of communication.
Suppose that Alice wants to tell Bob about a curve.
By Torelli's Theorem, she can do so by telling him the period matrix,
but this means transmitting $O(g^2)$ complex numbers in order to
describe something that depends on $O(g)$ parameters.
In other words, the period matrix is
{\sl sparse\/} in the sense of \cite{D}, and should therefore be compressible.  

The perspective that the period matrix is a compressible signal is
the central idea of the {\sl Information--Theoretic Schottky Problem\/}.
The attempt to apply ideas from Compressed Sensing \cite{D} to the
Schottky problem has led to many interesting experiments, conjectures,
and theorems \cite{W12}.  

The result described here is purely mathematical, rather than 
computational; however, it was  inspired by an attempt to implement ideas
in blind Compressed Sensing, as described in \cite{GE}.  

%%%%%%%%%%%%%%%%%%%%%%%%%%%%%%%%%%%%%%%%%%%%
\section{Plane Curves}

Shift the focus to a smooth plane curve whose affine
equation $f(x,y) = 0$ has degree $d > 4$.  Its genus is
$g = \frac{(d-1)(d-2)}{2}$, %g = {{(d-1)(d-2)}\over 2}$
and  its holomorphic differentials are given by 
$h(x,y) \frac{dx}{{\partial f}/{\partial y}},$
%{{dx}
%\over
%{
% {{\partial f}\over{\partial y}}
%}}$, 
where $h$ is
a so--called  adjoint polynomial of degree $d-3$.  Fix an order for the monomials
of degree $d-3$, {\it eg\/} the usual lexicographic order
$x^0y^0 < x^1y^0 < x^0y^1 < \cdots < x^0y^{d-3}$, thus
forming a proxy basis for $H^{(1,0)}(X)$.

%the 1--cycles on $X$ are more difficult to characterize, but 

The main result is

\begin{thm}\label{thm-main}
There is a set of 4 columns of the period matrix of a
smooth  plane curve that characterize the curve; in other words, if $X'$ is another
plane curve whose period matrix includes these four columns,
then $X$ and $X'$ are holomorphically equivalent.
\end{thm}

The four columns involved have $4g$ entries, so constitute a rather
small superset of ``moduli."  Thus, this is a significant loss--less compression
of the period matrix.

The number 4 seems rather arbitrary, but the condition that a curve have
a smooth planar representation is strong; one would not expect 
such a strong result from weaker hypotheses.

%%%%%%%%%%%%%%%%%%%%%%%%%%%%%%%%%%
\section{Period Matrices and Moduli}

The primary tool relating period matrices to moduli is the following
theorem of Rauch.  Let $K$ denote the canonical divisor on $X$.

\begin{thm}\label{thm-rauch} [Rauch]
Let $\{\zeta_1,\ldots \zeta_g\}$ be a normalized basis for $H^{(1,0)}(X)$
of a non-hyperelliptic Riemann surface $X$,
and suppose that $\{\zeta_i \zeta_j : (i,j) \in (I,J)\}$
form a basis for the quadratic differentials $H^0(X, 2K)$.
If another Riemann surface $X'$ has the same entries as $X$ in the
$(I,J)$ positions of its period matrix then $X$ and $X'$ are 
holomorphically equivalent.
\end{thm}

The proof, while not strictly relevant here, may be of interest for further
ITSP investigations.  It chooses the minimal member of the 
homotopy class of maps from the underlying surface of $X$ to
the underlying surface of $X'$ with respect to the Douglas--Dirichlet energy,
and proceeds by a delicate argument using infinitessimal quadratic
differentials to show that this map is holomorphic.

\bigskip

In principle, then, Alice can send Bob $3g-3$ entries of the period
matrix, and he can then verify that his period matrix is the same.
However, there is no canonical way to choose which $3g-3$ elements to
send, and there are many choices of $3g-3$ elements that
do not form moduli.  The point of Theorem A, then, is that in
the case of a smooth plane curve Alice can canonically choose a slightly
larger set of periods to send.

Returning to plane curves, the strategy is to choose a basis
for $H^0(2K)$ carefully; in the end, this will involve only
4 columns of the period matrix.

%%%%%%%%%%%%%%%%%%%%%%%%%%%
\section{Proof of the Theorem}

Recall the theorem of Noether (quoted in \cite{R}; also see \cite{ACGH}) that every quadratic
differential is a product of ordinary differentials.

To determine $(I,J)$, define a $g\times g$ matrix  $Q$ whose rows and columns are indexed
by the adjoint monomials.  In writing the matrix the factor 
$dx/\frac{\partial f}{\partial y}$ is omitted, and in considering
quadratic differentials one only need to look at products of the monomials.
In the case of degree $d=6$, the curve has genus
$g=10$ and, filling only the top row and leftmost column,
%
%$$Q = 
%\left [
%\matrix{
%\begin{matrix}
% 1 && x && y && x^2 && xy && y^2 && x^3 && x^2y && xy^2 && y^3\cr
% x &&    &&    &&        &&      && && && && && \cr
% y  && && &&  && && && && && &&\cr
% x^2 &&    &&    &&        &&      && && && && && \cr
 % xy &&    &&    &&  \ldots   &&      && && && && && \cr
%  y^2 &&    &&    &&        &&      && && && && && \cr
%  x^3 &&    &&    &&        &&      && && && && && \cr
%  x^2y &&    &&    &&        &&      && && && && && \cr
%  xy^2 &&    &&    &&        &&      && && && && && \cr
%  y^3 &&    &&    &&        &&      && && && && && \cr    
%  }
%\right ].
%\end{matrix}
%$$
\begin{displaymath}
Q = 
\left [
\begin{array}{cccccccccc}
 1 & x & y & x^2 & xy & y^2 & x^3 & x^2y & xy^2 & y^3 \\
 x &    &    &        &      & & & & & \\
 y  & & &  & & & & & &\\
 x^2 &    &    &        &      & & & & & \\
 xy &    &    &  \ldots   &      & & & & & \\
y^2 &    &    &        &      & & & & & \\
x^3 &    &    &        &      & & & & & \\
x^2y &    &    &        &      & & & & & \\
xy^2 &    &    &        &      & & & & & \\
y^3 &    &    &        &      & & & & & \\
\end{array} \right ]
\end{displaymath}

In this case, the columns beginning with 1, $x^3$, and $y^3$ must be 
included to form a basis for $H^0(X,2K)$.
More generally, the $x^{d-3}$ and $y^{d-3}$ columns must be included
in order to get all of the  monomials with $x$--degree (resp.~$y$--degree)
greater than $d-3$.  Note that these three columns contain
duplicate entries, for example, $x^{d-3}y^{d-3}$ appears in both
the $x^{d-3}$-- and $y^{d-3}$--columns.

The entries in these three columns do not constitute a basis, since they omit
the monomials of degree greater than $d-3$ but of $x$--
or $y$--degree $\le d-4$, but all of these monomials are in the
$x^2y^2$ column. To see this, let $x^ry^{d-2-r}$ be a monomial of
degree $d-2$; here $r \le d-4$.  This monomial factors
as $x^2y^2\cdot x^{r-2}y^{d-4-r}$, and  $x^{r-2}y^{d-4-r}$ is a monomial
in the first column.  Thus, $x^ry^{d-2-r}$ appears in the $x^2y^2$
column; the same applies to $x^{d-2-r}y^r$.

Similarly, monomials of degree $d-1$ can be written $x^ry^{d-1-r}$,
which factors as $x^2y^2 \cdot x^{r-2}y^{d-3-r}$.  Clearly $d - 3 - r \le
d-3$, so again such a monomial is a product of $x^2y^2$ and a monomial
from the first column.

The largest--degree monomial satisfying the conditions is
$x^{d-4}y^{d-4} = x^2y^2\cdot x^{d-6}y^{d-6}$. 

Thus, every ``missing" monomial appears in the $x^2y^2$ column.

Now consider the corresponding entries in the period matrix.  Since
the differentials of the chosen basis are not normalized, multiply
the right half of the period matrix by the inverse of the left half.  Each entry from 
$(I,J)$ in the normalized period matrix is a thus linear combination of entries from the
corresponding column in the matrix associated with $Q$.  But these entries still
correspond to a (superset) of a basis for $H^0(2K)$, and
thus by the Rauch Theorem determine the curve up to isomorphism.

%%%%%%%%%%%%%%%%%%%%%
\section{Complements}

$\bullet$ The four columns contain no more than $4g$ entries, which is
a substantial compression of the period matrix.

Even removing the duplicates, there are other relations among the
quadratic differentials and thus more relations between the periods.
This is easiest to see in degree $6$, where removing duplicate entries
leaves 28 positions, while the number of moduli is 27.  The missing
relation occurs in degree 6, and is, in fact, the equation of the 
curve.  In other words, some of the redundancy from 
the superset of periods used to determine the curve come from the
equation itself.

In higher degrees, many of the redundancies are in the ideal generated
by the equation.

$\bullet$ Since the columns of the normalized period matrix each correspond
to integrals over a cycle, it appears that only four of the generators of
the first singular homology group determine the whole topology of
the curve, but this is not the case because of the symmetry of the normalized
period matrix.

$\bullet$ It is neither true nor claimed that every set of four columns determines
the curve; using the Alice--Bob scenario, Alice, knowing that she
has a plane curve, chooses the columns used in the proof
and sends them to Bob.  Bob also has a curve, or perhaps a period
matrix, but may or may not know {\it a priori\/} whether his period
matrix is a plane curve, but if it contains the four columns then he has
determined that the curve that Alice sent is the one he has.  In this sense
the theorem provides more of a verification than an actual 
communication.

$\bullet$ D.~Litt points out that it may not be possible to transmit periods in
a finite message, although many complex numbers
do have compact descriptions ({\it eg\/} Gaussian rationals, surds).  In other
cases it may only be possible to transmit an approximation of the
periods.  If this is so, then Bob knows that his curve is close 
(in an analytic sense) to the
plane curve locus, which is already significant.

%
%%%%%%%%%%%%%%%%%%%%%%%%%%%%%%%%

\end{document}